\documentclass{amsart}
\usepackage{amssymb}
\usepackage{mathptm}

\usepackage[all]{xy}

\input xy
\xyoption{all}

\font\got=eufm10
\font\gots=eufm10 at 7pt
\numberwithin{equation}{section}
\newtheorem{lem}{Lemma}[section]
\newtheorem{corol}[lem]{Corollary}
\newtheorem{theor}[lem]{Theorem}
\newtheorem{prop}[lem]{Proposition}
\newtheorem{rem}[lem]{Remark}

\newtheorem{exm}[lem]{Example}
\newtheorem{examples}[lem]{Examples}

\def\natu{{\mathbb N}}
\def\entero{{\mathbb Z}}
\def\L{L_K(E)}
\def\M{\hbox{\got M}}
\def\Mp{\hbox{\gots M'}}
\def\F{\hbox{\got F}}

\def\soc#1{\mathop{\rm Soc}(#1)}
\def\a{\alpha}

\begin{document}
\title[Socle Theory for Leavitt path algebras of arbitrary graphs]{Socle Theory for Leavitt path algebras
of arbitrary graphs}

\author{G. Aranda Pino}
\address{Departamento de \'Algebra, Geometr\'{\i}a y Topolog\'{\i}a, Universidad de M\'alaga, 29071
M\'alaga, Spain.} \email{gonzalo@agt.cie.uma.es}

\author{D. Mart\'{\i}n Barquero}
\address{Departamento de Matem\'atica Aplicada, Universidad de M\'alaga, 29071
M\'alaga, Spain.}
\email{dmartin@uma.es}

\author{C. Mart\'{\i}n Gonz\'alez}
\address{Departamento de \'Algebra, Geometr\'{\i}a y Topolog\'{\i}a, Universidad de M\'alaga, 29071
M\'alaga, Spain.}
\email{candido@apncs.cie.uma.es}

\author{M. Siles Molina}
\address{Departamento de \'Algebra, Geometr\'{\i}a y Topolog\'{\i}a, Universidad de M\'alaga, 29071
M\'alaga, Spain.}
\email{msilesm@uma.es}


\begin{abstract}
The main aim of the paper is to give a socle theory for Leavitt path algebras of arbitrary graphs.
We use both the desingularization process and combinatorial methods
to study Morita invariant properties concerning the socle and to characterize it, respectively.
 Leavitt path algebras with nonzero socle are described as those which have line points, and it is 
 shown that the line points generate the socle of a Leavitt path algebra, extending so the results for row-finite
graphs in the previous paper \cite{AMMS} (but with different methods). A concrete description of
the socle of a Leavitt path algebra is obtained: it is a direct sum of matrix rings (of finite or infinite size) over the base field.

New proofs of the Graded Uniqueness and of the Cuntz-Krieger Uniqueness Theorems are given, shorthening
significantly the original ones. 
\end{abstract}

\maketitle

%%%%%%%%%%%%%%%%%%%%%%%%%%%%%%%%%%%%%%%%%%%%%%%%%%%%%%%%%%%%%%%%%%%%%%%%%%%%%%%%%%%%%%%%%%
%%%%%%%%%%%%%%%%%%%%%%%%%%%%%%%%%%%%%%%%%%%%%%%%%%%%%%%%%%%%%%%%%%%%%%%%%%%%%%%%%%%%%%%%%%
\section*{Introduction}

Leavitt path algebras of row-finite graphs have been recently introduced in \cite{AA1} and \cite{AMP}. They
have become a subject of significant interest, both for algebraists and for analysts working in C*-algebras.
The Cuntz-Krieger algebras $C^*(E)$ (the C*-algebra counterpart of these Leavitt path algebras) are described
in \cite{R}. The algebraic and analytic theories, while sharing some striking similarities, present some
remarkable differences, as was shown for instance in the ``Workshop on Graph Algebras'' held at the University
of M\'alaga (see \cite{Work}), and more deeply in the subsequent enlightening work by Tomforde
\cite{Tomforde}.

For a field $K$, the algebras $L_K(E)$ are natural generalizations of the algebras investigated by Leavitt in
\cite{Le}, and are a specific type of path $K$-algebras associated to a graph $E$ (modulo certain relations).
The family of algebras which can be realized as the Leavitt path algebra of a graph includes matrix rings
${\mathbb M}_n(K)$ for $n\in \mathbb{N}\cup \{\infty\}$ (where ${\mathbb M}_\infty(K)$ denotes matrices of
countable size with only a finite number of nonzero entries), the Toeplitz algebra, the Laurent polynomial
ring $K[x,x^{-1}]$, and the classical Leavitt algebras $L(1,n)$ for $n\ge 2$. Constructions such as direct
sums, direct limits and matrices over the previous examples can be also realized in this setting. But, in
addition to the fact that these structures indeed contain many well-known algebras, one of the main interests
in their study is the comfortable pictorial representations that their corresponding graphs provide.

The development of the theory of Leavitt path algebras (as well as that of their analytic sisters, the graph
C*-algebras) has had several different stages as far as questions of cardinality of the graphs are concerned.
At first, in the C*- case, only finite graphs (represented by matrices) were considered: Cuntz \cite{C} constructed and
investigated the C*-algebras ${\mathcal O}_n$ (nowadays called the Cuntz algebras), showing, among other
things, that each ${\mathcal O}_n$ is (algebraically) simple. Soon after the appearance of \cite{C}, Cuntz and
Krieger \cite{CK} described the significantly more general notion of the C*-algebra of a (finite) matrix $A$,
denoted ${\mathcal O}_A$. Among this class of C*-algebras one can find, for any finite graph $E$, the
Cuntz-Krieger algebra $C^*(E)$, defined originally in \cite{KPR}. The algebraic counterpart of these finite
Cuntz-Krieger algebras was considered in \cite{AGGP}.

The second step was to consider possibly infinite but countable row-finite graphs (that is, graphs with a
countable number of vertices and edges which satisfy that a vertex in the graph emits at most a finite number
of edges). This was first done in the analytic setting (see \cite{BPRS,R,RS} among others), while the seminal
results on Leavitt path algebras of row finite graphs appeared in \cite{AA1} and \cite{AMP}, so starting 
 a flurry of activity. In both situations the
classification of simple (\cite{AA1}) and purely infinite simple  (\cite{AA2}) structures was carried out in terms of properties of the
graph. In the analytic situation, the gauge invariant ideals were determined; in the algebraic one, the
graded ideals were described. In addition, several other ring properties were studied in the case of
row-finite Leavitt path algebras, such as being exchange \cite{APS}, finite dimensional \cite{AAS1},
noetherian \cite{AAS2}, semisimple \cite{AAPS} or having stable rank \cite{APS, AP}. It has been also shown that the Leavitt path algebra
$L_K(E)$ of a graph $E$ and the path algebra $KE$ associated to the same graph are closely enough: $L_K(E)$ is an algebra
of right quotients of $KE$ (see \cite{SilesAQLPA}).

Apart from the very recent paper \cite{Goodearl},  by K. R.  Goodearl, where he has introduced Leavitt path algebras of uncountable directed graphs,
the following breakthrough was to remove the hypothesis of row-finiteness in the underlying graphs.
Once more this was first done for the case of graph C*-algebras $C^*(E)$ (see for instance \cite{FLR,BHRS,DT})
and afterwards for Leavitt path algebras in \cite{AA3,Tomforde}. In both, the analytic and the algebraic
cases, very often the validity of the results for these not necessarily row-finite graphs requires coming
up with totally different proofs to those given for the row-finite case. Sometimes, helpful shortcuts such as
the desingularization process are at hand and sometimes the proofs must be reinvented. 
In this paper we face both situations.

Concretely, the aim of this article is to extend the theory of the socle of a Leavitt path algebra (considered
in \cite{AMMS} for row-finite graphs) to (countable) not necessarily row-finite graphs.
Specifically, we determine the structure of minimal left ideals of Leavitt path algebras of arbitrary graphs
and we scrutinize the nature of the socle of a Leavitt path algebra of an arbitrary graph in two different
ways: first, by singling out the set of vertices that generate the socle as a two-sided ideal (graph
description) and secondly, by unveiling the internal algebraic structure of it (algebraic description).

It is worth mentioning that these results on the socle were successfully applied for the row-finite case in
\cite{AAPS} in order to completely classify the semisimple/locally noetherian/locally artinian Leavitt path
algebras. Hence, this extension of \cite{AMMS} to arbitrary graphs could potentially led to achievements
analogous to \cite{AAPS} but for arbitrary graphs.

The article is organized as follows. The first section includes the basic definitions and examples that will
be used throughout. In addition, we describe several basic results and relations between the path algebra
and the Leavitt path algebra of an arbitrary graph.

In Section 2, a first approach to the study of the socle of a Leavitt path algebra of arbitrary graphs via the
desingularization process is made. We relate the ``line point" vertices of a graph (these are the vertices
that generate the socle as an ideal) to that of its desingularization. This allows us to establish, in
Corollary \ref{SocmorLeavitt}, further socle-related connections between the Leavitt path algebra of an
arbitrary graph $E$ and the Leavitt path algebra of its desingularized (row-finite) graph $F$. These are:
$L_K(E)$ has nonzero socle (respectively, coincides with its socle) if and only if $L_K(F)$ has nonzero socle
(respectively, coincides with its socle).

The action of the multiplication algebra is considered in Section 3. This is a valuable tool that allowed us
to shorten the lengthy proofs of the Graded Uniqueness Theorem and the Cuntz-Krieger Uniqueness Theorem for
Leavitt path algebras of arbitrary graphs given in \cite{Tomforde}. What is more, this tool allowed us to
weaken, in Theorem \ref{variante}, the set of hypotheses of the aforementioned Cuntz-Krieger Uniqueness
Theorem to a level that was useful for instance in \cite[Proof of Proposition 5.1]{AA3}.
These uniqueness theorems allow us to obtain that for a graph, the socle
of the Leavitt path algebra can be seen inside the socle of the corresponding graph C*-algebra. 

The study of minimal left ideals of Leavitt path algebras of arbitrary graphs is carried out in Section 4.
There are two phases for this: first, only principal left ideals generated by vertices are considered;
afterwards, principal left ideals generated by an arbitrary element of $L_K(E)$ are determined. One of the key
tools in \cite{AMMS} for this same enterprise in the row-finite case was the possibility to see the algebra
$L_K(E)$ as a certain direct limit of algebras associated to finite complete subgraphs. However, this
construction is no longer available in the general case of arbitrary graphs, so a completely different
approach, of a more combinatorial nature, is needed here.

Having paved the way, the natural subsequent and final step is taken in Section 5, where the socle of a
Leavitt path algebra of an arbitrary graph is determined. Thus, the graph description of the socle is given in
Theorem \ref{socleofLPA} and is the following: the socle of a Leavitt path algebra is the two-sided ideal
generated by the vertices of the graph whose trees (that is, the vertices which follow them in the graph) do
not contain cycles nor bifurcations (i.e., by line points). On the other hand, the algebraic description of the socle is given in
Theorem \ref{estructuradelzocalo}: the socle of a Leavitt path algebra is a direct sum of full matrix algebras
over the field $K$ of either finite or countable infinite size.

%%%%%%%%%%%%%%%%%%%%%%%%%%%%%%%%%%%%%%%%%%%%%%%%%%%%%%%%%%%%%%%%%%%%%%%%%%%%%%%%%%%%%%%%%%
%%%%%%%%%%%%%%%%%%%%%%%%%%%%%%%%%%%%%%%%%%%%%%%%%%%%%%%%%%%%%%%%%%%%%%%%%%%%%%%%%%%%%%%%%%
\section{Path algebras and Leavitt path algebras of arbitrary graphs}

A \emph{(directed) graph} $E=(E^0, E^1, r, s)$ consists of two countable sets $E^0$ and $E^1$, and maps $r, s
: E^1\to E^0$. The elements of $E^0$ are called \emph{vertices} and the elements of $E^1$ \emph{edges}. If for
a  vertex $v$, the set $s^{-1}(v)$ is finite, then the graph is called \emph{row-finite}. If $E^0$ is finite
then, by the row-finite hypothesis, $E^1$ must necessarily be finite as well; in this case we say simply that
$E$ is \emph{finite}. A vertex is called a \emph{sink} if it does not emit edges, and a \emph{source} if it
does not receive edges. A vertex $v$ such that $|s^{-1}(v)| = \infty$  is called an \emph{infinite emitter}.
Following \cite{Tomforde}, if $v$ is either a sink or an infinite emitter, we call it a \emph{singular
vertex}. If $v$ is not a singular vertex, we will say that it is a \emph{regular vertex}. A \emph{path} $\mu$
in a graph $E$ is a sequence of edges $\mu=e_1\dots e_n$ such that $r(e_i)=s(e_{i+1})$ for $i=1,\dots,n-1$. In
this case, $s(\mu):=s(e_1)$ is the \emph{source} of $\mu$, $r(\mu):=r(e_n)$ is the \emph{range} of $\mu$, and
$n$ is the \emph{length} of $\mu$, i.e, $l(\mu)=n$. We denote by $\mu^0$ the set of its vertices, that is:
$\mu^0=\{s(e_1),r(e_i):i=1,\dots,n\}$.

An edge $e$ is an {\it exit} for a path $\mu = e_1 \dots e_n$ if there exists $i$  such that $s(e)=s(e_i)$ and
$e\neq e_i$.  If $\mu$ is a path in $E$, and if $v=s(\mu)=r(\mu)$, then $\mu$ is called a \emph{closed path
based at $v$}. If $s(\mu)=r(\mu)$ and $s(e_i)\neq s(e_j)$ for every $i\neq j$, then $\mu$ is called a
\emph{cycle}.

We say that a graph $E$ satisfies \emph{Condition} (L) if every cycle in $E$ has an exit. For $n \geq 2$ we
define $E^n$ to be the set of paths of length $n$, and $E^\ast = \bigcup_{n\geq 0} E^n$ the set of all paths.

The set $T(v)=\{w\in E^0\mid v\ge w\}$ is the \emph{tree} of $v$, that is, the set of all the vertices in the
graph $E$ which follow $v$ ($v\ge w$ means that there is a path $\mu$ with $s(\mu)=v$ and $r(\mu)=w$). We will denote it by $T_E(v)$ when it is necessary to emphasize the dependence on
the graph $E$.

Now let $K$ be a field and let $KE$ denote the $K$-vector space which has as a basis the set of paths. It is
possible to define an algebra structure on $KE$ as follows: for any two paths $\mu=e_1\dots e_m, \nu=f_1\dots
f_n$, we define $\mu\nu$ as zero if $r(\mu)\neq s(\nu)$ and as the path $e_1\dots e_m f_1\dots f_n$ otherwise.
This $K$-algebra is called the \emph{path algebra of} $E$ \emph{over} $K$.

For a field $K$ and a  graph $E$, the \emph{Leavitt path $K$-algebra} $L_K(E)$ is defined as the $K$-algebra
generated by a set $\{v\mid v\in E^0\}$ of pairwise orthogonal idempotents, together with a set of variables
$\{e,e^*\mid e\in E^1\}$, which satisfy the following relations:

(1) $s(e)e=er(e)=e$ for all $e\in E^1$.

(2) $r(e)e^*=e^*s(e)=e^*$ for all $e\in E^1$.

(3) $e^*e'=\delta _{e,e'}r(e)$ for all $e,e'\in E^1$.

(4) $v=\sum _{\{ e\in E^1\mid s(e)=v \}}ee^*$ for every regular vertex $v\in E^0$.

The elements of $E^1$ are called \emph{(real) edges}, while for $e\in E^1$ we call $e^\ast$ a \emph{ghost
edge}.  The set $\{e^*\mid e\in E^1\}$ will be denoted by $(E^1)^*$.  We let $r(e^*)$ denote $s(e)$, and we
let $s(e^*)$ denote $r(e)$. If $\mu = e_1 \dots e_n$ is a path, then we denote by $\mu^*$ the element $e_n^*
\dots e_1^*$ of $L_K(E)$.

Formally speaking we can say that $L_K(E)$ is the quotient of the free associative algebra generated by
$E^0\cup E^1\cup (E^1)^*$ modulo the ideal induced by the identities (1)-(4). Hence the universal property of
the free associative algebra jointly with that of the quotient algebra can be used to construct homomorphisms
with domain $L_K(E)$.

We will recall some facts that will be used freely along the paper.

The Leavitt path $K$-algebra $L_K(E)$ is spanned as a $K$-vector space by $\{pq^* \mid p,q$ are paths in
$E\}$ (see \cite[Lemma 3.1]{Tomforde}). Moreover, $L_K(E)$ has a natural $\mathbb{Z}$-grading (see
\cite[Section 3.3]{Tomforde}): for each $n\in \mathbb{Z}$, the degree $n$ component $L_K(E)_n$ is spanned by
elements of the form $pq^*$ where $l(p)-l(q)=n$.

The set of \emph{homogeneous elements} is $\bigcup_{n\in {\mathbb Z}} L_K(E)_n$, and an element $x\in
L_K(E)_n$ is said to be $n$-\emph{homogeneous} or \emph{homogeneous of degree} $n$, denoted by $deg(x)=n$.

The $K$-linear extension of the assignment $pq^* \mapsto qp^*$ (for $p,q$ paths in $E$) yields an involution
on $L_K(E)$, which we denote simply as ${}^*$.  Clearly $(L_K(E)_n)^* = L_K(E)_{-n}$ for all $n\in {\mathbb
Z}$.

\begin{examples} {\rm By considering some basic configurations one can realize many algebras as
the Leavitt path algebra of some graph. Thus, for instance, the ring of Laurent polynomials $K[x,x^{-1}]$ is
the Leavitt path algebra of the graph

$$
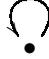
\xymatrix{{\bullet} \ar@(ur,ul)  }$$

Matrix algebras $M_n(K)$ can be achieved by considering a line graph with $n$ vertices and $n-1$ edges

$$\xymatrix{{\bullet} \ar [r]  & {\bullet} \ar [r]  & {\bullet} \ar@{.}[r] & {\bullet} \ar [r]  & {\bullet} }$$

Classical Leavitt algebras $L(1,n)$ for $n\geq 2$ are obtained as $L(R_n)$, where $R_n$ is the rose with $n$
petals graph

$$
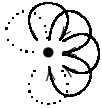
\xymatrix{{\bullet} \ar@(ur,dr)  \ar@(u,r)  \ar@(ul,ur)  \ar@{.} @(l,u) \ar@{.} @(dr,dl) \ar@(r,d) &  }$$
 \smallskip

Of course, combinations of the previous examples are possible. For example, the Leavitt path algebra of the
graph

$$
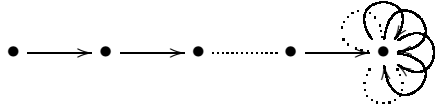
\xymatrix{{\bullet} \ar [r]  & {\bullet} \ar [r]  & {\bullet}
\ar@{.}[r] & {\bullet} \ar [r]  & {\bullet}
 \ar@(ur,dr)  \ar@(u,r)  \ar@(ul,ur)  \ar@{.} @(l,u) \ar@{.} @(dr,dl) \ar@(r,d) & }$$ \smallskip
 
\noindent
is $M_n(L(1,m))$, where $n$ denotes the number of vertices in the graph and $m$ denotes the number of loops.
In addition, the algebraic counterpart of the Toeplitz algebra $T$ is the Leavitt path algebra of the graph
$E$ having one loop and one exit
$$
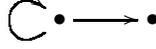
\xymatrix{{\bullet} \ar@(dl,ul) \ar[r] & {\bullet}  }$$
}
\end{examples}

There exists a natural inclusion of the path algebra $KE$ into the Leavitt path algebra $L_K(E)$ sending
vertices to vertices and edges to edges. We will use this monomorphism without any explicit mention to it.
Moreover, this natural monomorphism from the path algebra $KE$ into the Leavitt path algebra $L_K(E)$ is
graded, hence $KE$ is a $\mathbb{Z}$-graded subalgebra of $L_K(E)$.

We will revisit some basic results on the path algebra $KE$ and on the Leavitt path algebra $L_K(E)$ for an
arbitrary graph $E$. Given that the only difference between the definition of the Leavitt path algebra of a
row-finite graph and of an arbitrary graph is the non-existence of a CK2 relation ((4) in the definition) at
infinite emitters, it is perhaps not surprising that many of the results that hold for the row-finite case
still hold in this more general situation. In particular, by rereading the result in \cite[Lemma
1.1]{SilesAQLPA} we get that the following  still holds in this general situation.

\begin{lem}\label{indep} Let $E$ be an arbitrary graph. Any set of different paths is $K$-linearly independent.
\end{lem}

The same can be said about \cite[Proposition 2.1]{SilesAQLPA}:

\begin{prop}\label{semiprimeness}
Let $E$ be an arbitrary graph. Then the path algebra $KE$ is semiprime if and only if for every path $\mu$
there exists a path $\nu$ such that $s(\nu)=r(\mu)$ and $r(\nu)=s(\mu)$.
\end{prop}

Also, following the proof of \cite[Proposition 2.2]{SilesAQLPA} and applying \cite[Lemma 1.12]{ASgraded}
we have:

\begin{prop}\label{rgqa} For an arbitrary graph $E$ and any field $K$, the Leavitt path algebra $L_K(E)$
is an algebra of right quotients of the path algebra $KE$, equivalently, it is a $\mathbb{Z}$-graded algebra
of right quotients of $KE$.
\end{prop}

The following result was stated for row-finite graphs in \cite[Theorem 3.11]{AA1}. We include here a proof for
arbitrary graphs.

\begin{lem}\label{CicloSinSalidas} Let $E$ be an arbitrary graph. Let $v$ be a vertex in $E^0$ such that
there exists a cycle without exits $c$ based at $v$. Then:
$$vL_K(E)v= \left\{\sum_{i=-m}^n k_ic^i\ \vert \ k_i\in K;\ m, n \in \mathbb{N}\right\}\cong K[x, x^{-1}],$$
where $\cong$ denotes a graded isomorphism of $K$-algebras, and considering (by abuse of notation) $c^0=w$ and
$c^{-t}=({c^\ast})^t$, for any $t \geq 1$.
\end{lem}
\begin{proof}
First, it is easy to see that if $c=e_1 \ldots e_n $ is a cycle without exits based at $v$ and $u \in T(v)$,
then $s(f)=s(g)=u$, for $f,g \in E^1$, implies $f=g$. Moreover, if $r(h)=r(j)=w \in T(v)$, with $h,j \in E^1$, and
$s(h), \ s(j) \in T(v)$ then $h=j$. We have also that if $\mu \in E^*$ and $s(\mu )=u \in T(v)$ then there
exists $ k \in \natu^*, \ 1 \le k \le n$ verifying  $ \mu = e_k \mu'$ and $s(e_k)=u$.

Let $x \in vL_K(E)v$ be given by $x=\sum_{i=1}^p k_i \alpha_i\beta_i^*+\delta v$, with
$s(\alpha_i)=r(\beta_i^*)=s(\beta_i)=v$ and $\alpha_i, \ \beta_i \in E^*$. Consider $A=\{ \alpha \in E^*\colon
s(\alpha)=v \}$; we prove now that if $\alpha\in A$, $deg(\alpha)=mn+q, \ m,\ q \in \natu$ with $0\le q<n$, then
$\alpha=c^me_1 \dots e_q$. We proceed by induction on $deg(\alpha)$. If $deg(\alpha)=1$ and $s(\alpha)=s(e_1)$
then $\alpha=e_1$. Suppose now that the result holds for any $\beta\in A$ with $deg(\beta)\le sn+t$ and
consider any $\alpha\in A$, with $deg(\alpha)=sn+t+1$. We can write $\alpha=\alpha'f$ with $\alpha'\in A$, $f
\in E^1$ and $deg(\alpha')=sn+t$, so by the induction hypothesis $\alpha'=c^se_1 \dots e_t$.  Since
$s(f)=r(e_t)=s(e_{t+1})$ implies $f=e_{t+1}$, then $\alpha=\alpha'f=c^se_1 \dots e_{t+1}$.

We shall show that the elements $\alpha_i\beta_i^*$ are in the desired form, i.e.,  $c^d$ with $d \in \entero$.
Indeed, if $deg(\alpha_i)=deg(\beta_i)$ and $\alpha_i\beta_i^* \ne 0$,  we have $\alpha_i\beta_i^*=c^p
e_1\ldots e_k e_k^* \ldots e_1^*c^{-p}=v$ by (4). On the other hand $deg(\alpha_i)>deg(\beta_i)$ and
$\alpha_i\beta_i^* \ne 0$ imply $\alpha_i\beta_i^*=c^{d+q} e_1\ldots e_k e_k^* \ldots e_1^*c^{-q}=c^d, \ d \in
\natu^*$. In a similar way, from $deg(\alpha_i)<deg(\beta_i)$ and $\alpha_i\beta_i^* \ne 0$ it follows that
$\alpha_i\beta_i^*=c^{q} e_1\ldots e_k e_k^* \ldots e_1^*c^{-q-d}=c^{-d}, \ d \in \natu^*$. Define
$\varphi\colon K[x,x^{-1}] \to L_K(E)$ by $\varphi(1)=v, \ \varphi(x)=c$ and $\varphi(x^{-1})=c^*$. It is a
straightforward routine to check that $\varphi$ is a graded monomorphism with image $vL_K(E)v$, so that $vL_K(E)v$ is
graded isomorphic to $K[x,x^{-1}]$ as a graded $K$-algebra.
\end{proof}

%%%%%%%%%%%%%%%%%%%%%%%%%%%%%%%%%%%%%%%%%%%%%%%%%%%%%%%%%%%%%%%%%%%%%%%%%%%%%%%%%%%%%%%%%%
%%%%%%%%%%%%%%%%%%%%%%%%%%%%%%%%%%%%%%%%%%%%%%%%%%%%%%%%%%%%%%%%%%%%%%%%%%%%%%%%%%%%%%%%%%
\section{Desingularization, Morita equivalence and socle}
Given an arbitrary graph $E$, one can associate a row-finite graph $F$, called a desingula\-rization of $E$,
such that $L_K(E)$ is Morita equivalent to $L_K(F)$ as rings with local units \cite[Theorem 5.2]{AA3}. The
process of building the graph $F$ out of $E$ is described in \cite{DT, AA3} and it essentially consists, as
the name suggests, on conveniently removing the singular vertices of $E$. We briefly recall the process here
for the reader's convenience:

If $v_0$ is a sink in $E$, then by \emph{adding a tail at $v_0$} we mean attaching a graph of the form
$$\xymatrix{ {\bullet}^{v_0} \ar[r] & {\bullet}^{v_1} \ar[r] & {\bullet}^{v_2} \ar[r] &
 {\bullet}^{v_3} \ar@{.>}[r] & }$$
to $E$ at $v_0$. If $v_0$ is an infinite emitter  in $E$, then by \emph{adding a tail at $v_0$} we mean
performing the following process: we first list the edges $e_1, e_2, e_3, \ldots$ of $s^{-1}(v_0)$, then we
add a tail to $E$ at $v_0$ of the following form
$$\xymatrix{ {\bullet}^{v_0} \ar[r]^{f_1} & {\bullet}^{v_1} \ar[r]^{f_2} & {\bullet}^{v_2} \ar[r]^{f_3} &
{\bullet}^{v_3} \ar@{.>}[r] & }$$ We remove the edges in $s^{-1}(v_0)$, and for every $e_j \in s^{-1}(v_0)$
we draw an edge $g_j$ from $v_{j-1}$ to $r(e_j)$.

If $E$ is a directed graph, then a \emph{desingularization} of $E$ is a graph $F$ formed by adding a tail to
every sink and every infinite emitter of $E$ in the fashion above. Several basic examples of desingularized
graphs are found in \cite[Examples 5.1, 5.2 and 5.3]{AA3}.

\par

When extending results to Leavitt path algebras of arbitrary graphs two main philosophies have been followed.
The obvious one consists on just reproving the results for arbitrary graphs with some ad hoc methods, while
the second approach uses the aforementioned desingularization construction which allows us to transfer, via a
Morita equivalence, results from the arbitrary graph setting to the row-finite situation.

In this section, we will obtain information about the socle of an arbitrary Leavitt path algebra $L_K(E)$ out of
the information provided by the socle of the Leavitt path algebra of its desingularization $L_K(F)$.

A vertex $v$ in $E^0$ is a \emph{bifurcation} (or \emph{there is a bifurcation at} $v$) if $s^{-1}(v)$ has at
least two elements. A vertex $u$ in $E^0$ will be called a \emph{line point} if there are neither bifurcations
nor cycles at any vertex $w\in T(u)$. We will denote by $P_l(E)$ the set of all line points in $E^0$.

As we will show, the set of line points plays a crucial role in the description of the socle of a Leavitt path
algebra of an arbitrary graph $E$. The next results analyze the relation between the line points of $E$ and those
of the desingularized graph $F$.

\begin{prop}\label{desingpointline} Let $E$ be an arbitrary graph and $F$ any desingularization of $E$. Then
\begin{enumerate}
\item $P_l(E)=P_l(F)\cap E^0$.
\item $P_l(E)\neq \emptyset$ if and only if $P_l(F)\neq \emptyset$.
\end{enumerate}
\end{prop}
\begin{proof}
(1). Suppose that $v\in P_l(E)$. Then $T_E(v)$ does not contain bifurcations nor cycles in $E$; in particular,
it does not contain infinite emitters in $E$. Therefore, no edges are added at any vertex of $T_E(v)$ in the
desingularization process unless $T_E(v)$ contains a (necessary unique) sink $w$, in which case an infinite
tail of the form $$\xymatrix{ {\bullet}^{w} \ar[r] & {\bullet} \ar[r] & {\bullet} \ar[r] & {\bullet}
\ar@{.>}[r] & }$$ has been attached at $w$. In other words, $T_F(v)$ does not contain bifurcations nor cycles
in $F$ either, that is, $v\in P_l(F)$.

To see the converse containment, take $v\in P_l(F)\cap E^0$. Then $T_F(v)$ does not contain bifurcations nor
cycles in $F$. Note that, by construction, neither vertex of $T_F(v)$ is a sink nor an infinite emitter. All
this shows that there exists a countable family of vertices $\{v_i\}_{i=0}^\infty$ and edges
$\{e_i\}_{i=0}^\infty$ of $F$ such that $v_0=v$, $s_F^{-1}(v_i)=\{e_i\}$ for all $i$ and
$T_F(v)=\{v_i\}_{i=0}^\infty$.

Since $v\in E^0$ we have two situations. First, if every $v_i$ was already in $E$, then from the way the graph
$F$ is constructed we conclude that $T_E(v)=\{v_i\}_{i=0}^\infty$, and $s_E^{-1}(v_i)=\{e_i\}$ for all $i$.
Otherwise, there exists $j\geq 0$, such that $v_j$ is a sink in $E^0$, $T_E(v)=\{v_i\}_{i=0}^j$, and
$s_E^{-1}(v_i)=\{e_i\}$ for all $i\leq j$. In both cases $v\in P_l(F)$ implies that $v\in P_l(E)$.

\medskip

(2). If $P_l(E)\neq \emptyset$, then $P_l(F)\neq \emptyset$ by (1). Suppose now that $v\in P_l(F)$. Again, if
$v\in E^0$, then $v\in P_l(E)$ by (1). Otherwise, if $v$ is a vertex which was not originally if $E$, then it
cannot be a vertex in $\{v_i\}_{i\geq 1}$ of any new infinite tail of the form
$$
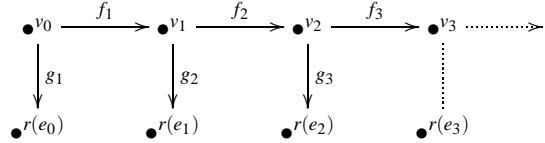
\xymatrix{ {\bullet}^{v_0} \ar[r]^{f_1} \ar[d]^{g_1} &
{\bullet}^{v_1} \ar[r]^{f_2} \ar[d]^{g_2} & {\bullet}^{v_2} \ar[r]^{f_3} \ar[d]^{g_3} & {\bullet}^{v_3} \ar@{.>}[r] \ar@{.}[d] &  \\
{\bullet}^{r(e_0)}  & {\bullet}^{r(e_1)} & {\bullet}^{r(e_2)} & {\bullet}^{r(e_3)} & }$$ as the trees of all
vertices $\{v_i\}_{i\geq 0}$ of these configurations necessarily have bifurcations. Therefore, $v$ is a vertex
of an infinite tail of $F$ which was introduced at a sink $z$ in $E$, but then clearly $z\in P_l(E)$. So that
$P_l(E)\neq \emptyset$ in this case too.
\end{proof}

One of the main results of \cite{AA3} was \cite[Theorem 5.2]{AA3}. There, the authors proved that if $E$ is an
arbitrary graph, then $L_K(E)$ is Morita equivalent to $L_K(F)$ for any desingularization $F$ of $E$. We are
going to exploit that fact in this section. First, we recall the notion of Morita equivalence for idempotent
rings (a ring $R$ is said to be \emph{idempotent} if $R^2=R$). Note that since Leavitt path algebras have
local units, they are idempotent rings.

Let $R$ and $S$ be two  rings, $_R N_S$ and $_S M_R$ two bimodules and $(-, -): N \times M \to R$, $[  -, -]:
M \times N \to S$ two maps. Then the following conditions are equivalent:
\begin{enumerate}
\item[(i)] $\left(\begin{matrix}
 R & N \cr M & S\end{matrix}\right)$
  is a ring with componentwise sum and product given by:
           $$ \left(\begin{matrix}r_1 & n_1 \cr m_1 & s_1\end{matrix}\right)
          \left(\begin{matrix} r_2 & n_2 \cr m_2 & s_2\end{matrix}\right) =
         \left(\begin{matrix}
         r_1r_2 +(n_1, m_2) & r_1n_2 +n_1s_2 \cr
          m_1r_2 +s_1m_2 & [  m_1, n_2] +s_1s_2
          \end{matrix}\right)$$
\item[(ii)] $[ -, -]$ is $S$-bilinear and $R$-balanced,  $( -, -)$ is $R$-bilinear and $S$-balanced
and the following associativity conditions hold:
$$(n, m)n^{\prime} = n [m, n^{\prime}] \quad \hbox{and} \quad
            [m, n] m^{\prime} = m (n, m^{\prime}).$$
            $[ -, -]$ being $S$-bilinear and $R$-balanced and  $ ( -, -)$ being $R$-bilinear and
           $S$-balanced is equivalent to having bimodule maps
        $\varphi : N \otimes_S M \to R$ and  $\psi : M \otimes_R N \to  S$, given by
       $$ \varphi (n \otimes m) = (n, m) \quad \hbox{and} \quad \psi (m \otimes n) = [m, n]$$
        so that the associativity conditions above read
       $$  \varphi (n \otimes m) n^\prime= n \psi (m \otimes n^\prime) \quad \hbox{and} \quad
       \psi (m \otimes n) m^\prime = m \varphi (n\otimes m^\prime).$$
\end{enumerate}
A \emph{Morita context} is a sextuple $(R, S, N, M, \varphi, \psi)$ satisfying the conditions given above. The
associated ring is called the \emph{Morita ring of the context}. By abuse  of notation we will write $(R, S,
N, M)$ instead  of $(R, S, N, M, \varphi, \psi)$ and  will suppose $R$, $S$, $N$, $M$ contained in the Morita
ring associated to the context. The Morita context will  be called \emph{surjective} if the maps $\varphi$ and
$\psi$ are both surjective.

In classical Morita theory, it is shown that two rings with identity $R$ and $S$ are Morita equivalent (i.e.,
$R$-mod and $S$-mod are equivalent categories) if and only if there exists a surjective Morita context $(R, S,
N, M, \varphi, \psi)$. The approach to  Morita theory for rings without identity by means of Morita contexts
appears in a number of papers (see \cite{GarciSimon} and the references therein) in which many consequences
are obtained from the existence of a Morita context for two rings $R$ and $S$.

For an idempotent ring $R$ we denote by $R-$Mod the full subcategory of the category of all left $R$-modules
whose objects are the ``unital" nondegenerate modules. Here, a left $R$-module $M$ is said to be \emph{unital}
if $M=RM$, and $M$ is said to be \emph{nondegenerate} if, for $m\in M$, $Rm=0$ implies $m=0$. Note that, if
$R$ has an identity, then $R-$Mod is the usual category of left $R-$modules $R$-mod.

It is shown in \cite[Theorem]{Ky} that, if $R$ and $S$ are arbitrary rings having a surjective Morita context,
then the categories $R-$Mod and $S-$Mod are equivalent. It is proved in \cite[Proposition 2.3]{GarciSimon}
that the converse implication holds for idempotent rings.

Given two idempotent rings $R$ and $S$, we will say that they are \emph{Morita equivalent} if the respective
full subcategories of unital nondegenerate modules over $R$ and $S$ are equivalent.

The following result can be found in \cite{GarciSimon} (see Proposition 2.5 and Theorem 2.7).

\medskip

\begin{theor}Let $ R $ and $ S $ be two idempotent rings. Then the categories $R-$Mod and
$S-$Mod are equivalent if and only if there exists a surjective Morita context $ (R, S, M, N) $.
\end{theor}
\medskip

The socles of Morita equivalent semiprime idempotent rings are closely related as we will see next. The proofs
of the following results are largely based on the concept of local algebra at an element that we proceed to
introduce.

For a ring $R$ and an element $x\in R$, the \textit{local ring} of $R$ \textit{at} $x$ (denoted $R_x$) is
defined to be the ring $xRx$, with the sum inherited from $R$, and product given by $xax\cdot xbx= xaxbx$.
The use of local rings at elements allows to overcome the lack of a unit element in the original ring, and to translate
problems from a non-unital context to the unital one.
See \cite{GS} for an equivalent definition and information about the exchange of properties between a ring
and its local rings at elements. In particular, if $e$ is an idempotent in the ring $R$, then the local ring
of $R$ at $e$ is just the corner $eRe$.

If $R$ is a semiprime ring, then the sum of all its minimal left ideals coincides with the sum of all its
minimal right ideals. This sum is called the socle of $R$ and will be denoted by $\soc{R}$. When the ring has
no minimal one-sided ideals, it is said that $R$ has zero socle.

As  the next lemma shows, taking local rings at elements and considering the socle are commuting operations.

\begin{lem}\label{zocalocal} For a semiprime ring $R$ and an element $x \in R$, we have
$\soc{R_x}= (\soc{R})_x$.
\end{lem}
\begin{proof} Note that for every element $x\in R$, being $R$ semiprime implies $R_x $ is semiprime too
(apply \cite[Proposition 2.1 (i)]{GS}), hence it has sense to consider the socle of the local ring at the
element.

Show first  $\soc{R_x} \subseteq (\soc{R})_x$. If $xax\in \soc{R_x}$, by \cite[Proposition 2.1 (v)]{GS}, $xax
\in \soc{R}$. As the socle is a von Neumann regular ring, there exists $y\in R$ such that
$xax=(xax)y(xax)=(xax)y(xax)y(xax)$; use that the socle is an ideal of the ring and that  $xax$ is in the
socle of $R$ to obtain $axyxaxyxa\in \soc{R}$, so that $xax= x (axyxaxyxa) x \in x (\soc{R}) x$.

For the converse, consider $xax\in x(\soc{R})x$, with $a \in \soc{R}$. Apply again that the socle is an ideal
to deduce that $xax$ is in $\soc{R}$. By \cite[Proposition 2.1 (v)]{GS} $xax \in \soc{R_x}$, as wanted.
\end{proof}

\begin{theor}\label{Socmor} Let $R$ and $S$ be two Morita equivalent semiprime idempotent rings. Then:
\begin{enumerate}
\item $R$ has nonzero socle if and only if $S$ has nonzero socle.
\item $R=\soc{R}$ if and only if $S=\soc{S}$.
\end{enumerate}
\end{theor}
\begin{proof} We start by setting several notation and results that will be used to prove the statements.

Denote by $A$ the Morita ring associated to a surjective Morita context $ (R, S, N, M) $, and identify $R$,
$S$, $N$ and $M$, in the natural way, with subsets of $A$.

(i) Use Lemma \ref{zocalocal} to settle that for any $x\in R$ we have: $ (\soc{R})_x=\soc{R_x}=\soc{A_x}=
 (\soc{A})_x$.

(ii) $\soc{R} = \soc{A} \cap R $ and $\soc{S} = \soc{A} \cap S $. This follows from (i) together with the fact
that $R_x=A_x$ for every $x\in R$ and the fact that an element is in the socle of a ring if and only if the
local ring at the element is an artinian ring (see \cite[Proposition 1.2 (v)]{GS}).

(1). Take a nonzero element $x$ in $\soc{R}$.  By (ii), $x\in \soc{A}$, and as the socle is an ideal, $MxN$,
which is contained in $S$, is in the socle of $A$ too.   We claim that $MxN$ is nonzero because otherwise
$0=NMxN=RxR$, a contradiction since every element in the socle is von Neumann regular and $x$ is nonzero.
Therefore we have $0\neq MxN\subseteq \soc{A}\cap S\subseteq \soc{S}$.

It can be proved, in an analogous way, that $\soc{S} \neq 0$ implies $\soc{R}\neq 0$.

(2). If $R$ coincides with its socle, $R\subseteq \soc{A}$. Then $S=MN=MNMN=MRN\subseteq M (\soc{A})
N\subseteq  \soc{A}\cap S=  \soc{S}$.

\end{proof}

The results above can be readily adapted to the Leavitt path algebra setting.

\begin{corol}\label{SocmorLeavitt} Let $E$ be an arbitrary graph and $F$ any desingularization of $E$. Then
\begin{enumerate}
\item $\soc{L_K(E)}\neq 0$ if and only if $\soc{L_K(F)}\neq 0$.
\item $L_K(E)$ coincides with its socle if and only if $L_K(F)$ coincides with its socle.
\end{enumerate}
\end{corol}
\begin{proof} Use first \cite[Theorem 5.2]{AA3} to get that $L_K(E)$ is Morita equivalent to $L_K(F)$. Now the
proof is a straightforward consequence of Theorem \ref{Socmor} and the fact that Leavitt path algebras for
arbitrary graphs are rings with local units (hence, idempotent rings), and also semiprime \cite[Proposition
6.1]{AA3}.
\end{proof}

The following result is a generalization of \cite[Corollary 4.3]{AMMS} for arbitrary graphs.

\begin{corol} Let $E$ be an arbitrary graph, then $L_K(E)$ has nonzero socle if and only if
$P_l(E)\neq \emptyset$.
\end{corol}
\begin{proof}
Consider any desingularization $F$ of $E$. Apply Corollary \ref{SocmorLeavitt} (1) to obtain that
$\soc{L_K(E)}\neq 0$ if and only if $\soc{L_K(F)}\neq 0$. By the row-finite case proved in \cite[Corollary
4.3]{AMMS} we know that $\soc{L_K(F)}\neq 0$ if and only if $P_l(F)\neq \emptyset$. Finally, use Proposition
\ref{desingpointline} (2) to get the result.
\end{proof}

One of the main aims of the paper is the complete determination of the socle of a Leavitt path algebra of an
arbitrary graph as the ideal generated by its set of line points. Unfortunately, this description of
the socle is unreachable via Morita equivalence (that is, by using a desingularization process).
Among other things, because two Morita equivalent idempotent rings can have socles of different size. For example:

\begin{exm}{\rm Let $R:=RCFM(K)$ be the ring of infinite matrices with entries in a field $K$, and finite rows and columns. Consider in $R$ the idempotent $e_{11}$ (defined as the matrix having 1 in place $(1,1)$ and zero elsewhere), and denote $f:=1-e\in  R$. Then $(eRe, fRf, fRe, eRf)$ is
a surjective Morita context for the two idempotent rings $eRe$ and $fRf$, and while $eRe$ is finite dimensional (in fact, it is isomorphic to the base field), $fRf$ is not.}
\end{exm}

In the upcoming sections of the paper we will use specifically-adapted methods in order
to achieve our main goal: the description of the socle.

%%%%%%%%%%%%%%%%%%%%%%%%%%%%%%%%%%%%%%%%%%%%%%%%%%%%%%%%%%%%%%%%%%%%%%%%%%%%%%%%%%%%%%%%%%
%%%%%%%%%%%%%%%%%%%%%%%%%%%%%%%%%%%%%%%%%%%%%%%%%%%%%%%%%%%%%%%%%%%%%%%%%%%%%%%%%%%%%%%%%%
\section{Action of the multiplication algebra of $\L$}

The action of the multiplication algebra of a Leavitt path algebra on the algebra itself has proved to be a
powerful tool which has allowed to shorten the lengthy proofs of some results in this theory. Recall that for
a not necessarily associative $K$-algebra $A$, and fixed $x,y\in A$, the left and right \emph{multiplication
operators} $L_x,R_y\colon A\to A$ are defined by $L_x(y):=xy$ and $R_y(x):=xy$. Denoting by $\hbox{End}_K(A)$
the $K$-algebra of $K$-linear maps $f\colon A\to A$, the \emph{multiplication algebra} of $A$ (denoted
$\M(A)$) is the subalgebra of $\hbox{End}_K(A)$ generated by the unit and all left and right multiplication
operators $L_a,R_a:A\to A$. There is a natural action of $\M(A)$ on $A$ such that $A$ is an $\M(A)$-module
whose submodules are just the ideals of $A$. This is given by $\M(A)\times A\longrightarrow A$, where $f\cdot
a:=f(a)$ for any $(f,a)\in\M(A)\times A$. Given $x,y\in A$ we shall say that $x$ \emph{is linked to} $y$ if
there is some $f\in\M(A)$ such that $y=f(x)$. This fact will be denoted by $x\vdash y$.

The result that follows was proved in \cite[Proposition 3.1]{AMMS} for row-finite graphs. It states that any
nonzero element in a Leavitt path algebra is linked to either a vertex or to a nonzero polynomial in a cycle
with no exits. So it gives a full account of the action of $\M(\L)$ on $\L$. This result proved to be very
powerful as the main ingredient to show that the socle of a Leavitt path algebra of a row-finite graph is the
ideal generated by the line points. The same proof given there can be used in the case of not necessarily
row-finite graphs.

\begin{prop}\label{zuru} Let $E$ be an arbitrary graph. Then, for every nonzero element $x\in L_K(E)$,
there exist $\mu_1,\dots,\mu_r,\nu_1,\dots,\nu_s\in E^0\cup E^1\cup (E^1)^*$ such that:
\begin{enumerate}
\item $\mu_1\dots\mu_rx\nu_1\dots\nu_s$ is a nonzero element in $Kv$, for some $v\in E^0$, or
\item there exist a vertex $w\in E^0$ and a cycle without exits $c$ based at $w$ such that
$\mu_1\dots\mu_rx\nu_1\dots\nu_s$ is a nonzero element in $wL_K(E)w$.
\end{enumerate}

Both cases are not mutually exclusive.
\end{prop}

\begin{corol}\label{corzuru}
For any nonzero $x\in\L$ we have $x\vdash v$ for some $v\in E^0$ or $x\vdash p(c,c^*)$ where $c$ is a cycle
with no exits and $p$ a nonzero polynomial in $c$ and $c^*$.
\end{corol}
\begin{proof}
Use Lemma \ref{CicloSinSalidas} together with Proposition \ref{zuru}.
\end{proof}

For any $K$-algebra $A$ the $\M(A)$-submodules of $A$ are just the ideals of $A$ and the cyclic
$\M(A)$-submodules of $A$ are the ideals generated by one element (\emph{principal ideals} in the sequel). So
the previous corollary states that the nonzero principal ideals contain either vertices or nonzero elements of
the form $p(c,c^*)$. Therefore, for graphs in which every cycle has an exit, each nonzero ideal contains a
vertex. Now, \cite[Corollary 6.10]{Tomforde} can be obtained immediately from item (ii) in the following
result:

\begin{corol}\label{VerticeEnTodoIdeal} Let $E$ be an arbitrary graph.
\begin{enumerate}
\item[(1)] Every ${\mathbb Z}$-graded nonzero ideal of $L_K(E)$ contains a vertex.
\item[(2)] Suppose that $E$ satisfies Condition {\rm(L)}. Then every nonzero ideal of $L_K(E)$ contains a
vertex.
\end{enumerate}
\end{corol}
\begin{proof}
The second assertion has been proved above. So assume that $I$ is a graded ideal of $\L$ which contains no
vertices. Let $0\neq x\in I$ and use Corollary \ref{corzuru} to find elements $y,z\in L_K(E)$ such that
$yxz=\sum_{i=-m}^n k_i c^i \neq0$. But $I$ being a graded ideal implies that every summand is in $I$. In
particular, for $t\in \{-m, \dots, n\}$ such that $k_tc^t \neq 0$ we have $0\neq (k_t)^{-1}c^{-t}k_t c^t =
w\in I$, which is absurd.
\end{proof}

It was shown in \cite[Proposition 6.1]{AA3} that $L_K(E)$ is a semiprime algebra for an arbitrary graph. 
The proof required the use of the desingularization process. Here, we can give
an element-wise proof by using Proposition \ref{zuru}.

\begin{prop} \label{nodegenerada} Let $E$ be an arbitrary graph. Then $L_K(E)$ is semiprime.
\end{prop}
\begin{proof}
Take a nonzero ideal $I$ such that $I^2=0$. If $I$ contains a vertex we are
done. On the contrary there is a nonzero element $p(c,c^*)\in I$ by Corollary \ref{corzuru}. If we consider
the (nonzero) coefficient of maximum degree in $c$ and write $p(c,c^*)^2=0$ we immediately see that this
scalar must be zero, a contradiction.
\end{proof}

To illustrate how powerful is Proposition \ref{zuru} we can see how it reduces considerably in
length the proofs given in \cite{Tomforde} of the so-called Uniqueness Theorems. These are \cite[Theorem
4.6]{Tomforde} (Graded Uniqueness Theorem) and \cite[Theorem 6.8]{Tomforde} (Cuntz-Krieger Uniqueness
Theorem).

\begin{theor} \label{uniqueness}
Let $E$ be an arbitrary graph, and let $L_K(E)$ be the associated Leavitt path algebra.
\begin{enumerate}
\item[(1)] Graded Uniqueness Theorem.

If $A$ is a $\mathbb{Z}$-graded ring and $\pi:L_K(E)\to A$ is a graded ring homomorphism with $\pi(v)\neq 0$
for every vertex $v\in E^0$, then $\pi$ is injective.

\item[(2)] Cuntz-Krieger Uniqueness Theorem.

Suppose that $E$ satisfies Condition {\rm(L)}. If $\pi: L_K(E)\to A$ is a ring homomorphism with $\pi(v) \neq
0$, for every vertex $v\in E^0$, then $\pi$ is injective.
\end{enumerate}
\end{theor}

\begin{proof} In both cases, the kernel of the ring homomorphism $\pi$ is an algebra ideal
(a graded ideal in the first one). By Corollary \ref{VerticeEnTodoIdeal},  $\textrm{Ker} (\pi)$ must be zero
because otherwise it would contain a vertex (apply (i) in the corollary to (1) and (ii) to the other case),
which is not possible by the hypotheses.
\end{proof}

In \cite{Tomforde}, Tomforde used the previous theorems in order to prove that, for the field of complex
numbers ${\mathbb C}$, the Leavitt path algebra $L_{\mathbb C}(E)$ could be embedded in the graph C*-algebra
$C^*(E)$ via a homomorphism $\phi:L_{\mathbb C}(E)\to C^*(E)$ sending the generators of $L_{\mathbb C}(E)$ to
the generators of $C^*(E)$. As he noted, such a homomorphism is well-defined by the universal property of
$L_{\mathbb C}(E)$ and is injective, as can be shown precisely by  applying the Graded Uniqueness Theorem. Here, we can
use this embedding $\phi$ to get that the socle of the Leavitt path algebra of an arbitrary graph $E$ is
always contained in the socle of the graph C*-algebra of $E$ (but may not be equal), as is shown in the next
result. The counterexample contained in the following proposition was communicated to the authors by Pere Ara.

\begin{prop}\label{contenidos} Let $E$ be an arbitrary graph. Then
$\soc{L_{\mathbb C}(E)}\subseteq \soc{C^*(E)}$. Moreover, there exists a row-finite graph $E$ such that the
inclusion is proper.
\end{prop}
\begin{proof}
As explained in the previous paragraph, we can use Theorem \ref{uniqueness} (1) and the ideas of \cite[Proof
of Theorem 7.3]{Tomforde} to obtain that $L_{\mathbb C}(E)$ is isomorphic to a dense *-subalgebra ${\mathcal
A}$ of $C^*(E)$. 

Consider a minimal idempotent $e\in {\mathcal A}$. We will show that $e$ is a minimal idempotent in $C^*(E)$
as well. It suffices to show that $eC^*(E)e$ is a division ring. Take a nonzero element $x\in eC^*(E)e$.
Because ${\mathcal A}$ is a dense *-subalgebra of $C^*(E)$, there exists a sequence
$\{x_n\}_{n=1}^\infty\subseteq {\mathcal A}$ such that $x=\lim (ex_ne)$. Suppose that $n$ is such that the
element $ex_ne\in e{\mathcal A}e$ is nonzero. Since $e{\mathcal A}e$ is a division ring, there exists $y_n\in
e{\mathcal A}e$ such that $ex_ne y_n=y_n ex_ne=1|_{e{\mathcal A}e}=e$. As $x\neq 0$, there exists $m$ such
that $ex_ne\neq 0$ for every $n\geq m$. Define $y_n=0$ for every $n<m$ and $y=\lim y_n$. Then $xy=\lim
(ex_ne)y_n=\lim e=e=1|_{eC^*(E)e}$ and analogously $yx=1|_{eC^*(E)e}$. This proves our claim.

Recall that $\soc{R}$ is the two-sided ideal generated by the minimal idempotents of a ring $R$. Denote by
${\mathcal I}$ the set of the minimal idempotents in $\mathcal A$, and by$\mathcal I^\ast$ the set of minimal idempotents
in $C^*(E)$. Then:
$$\soc{\mathcal A}  =  \sum_{e \in {\mathcal I}}
{\mathcal A}e{\mathcal A} \subseteq \sum_{e	\in {\mathcal I}} C^*(E)eC^*(E)
\subseteq  \sum_{e\in {\mathcal I^\ast}} C^*(E)eC^*(E) =
\soc{C^*(E)}.$$

To show that the inclusion might be proper we consider the row-finite graph $E$ given by
$$\xymatrix{ {\bullet} \ar@(dl,ul) \ar[r] & {\bullet}}$$ 

In \cite{SilesAQLPA}, Siles Molina showed that the
Leavitt path algebra of this graph is the algebraic Toeplitz algebra $T={\mathbb C}\langle x,y\ |\
xy=1\rangle$, for which it is known that $\soc{T}={\mathbb M}_\infty({\mathbb C})$. However, the completion of
$T$ is the analytic Toeplitz algebra, whose socle is the algebra of finite rank operators, which strictly
contains ${\mathbb M}_\infty({\mathbb C})$ (matrices of countable size with only a finite number of nonzero entries).
\end{proof}

In spite of the power of the aforementioned Uniqueness Theorems, one may encounter different situations in
which neither set of hypotheses in Theorem \ref{uniqueness} are satisfied. This happens, for instance, in the
proof of \cite[Proposition 5.1]{AA3}. Here, the authors showed that for any arbitrary graph $E$, and $F$ any
of its desingularizations, the Leavitt path algebra $L_K(E)$ is isomorphic to a subalgebra of $L_K(F)$. In
proving this result, they built a homomorphism $\phi$ from $L_K(E)$ to $L_K(F)$ and needed to show its
injectivity. In this situation, neither $E$ satisfied Condition (L) nor $\phi$ was a graded homomorphism
(because the desingularization process might enlarge some paths but not all of them), so the Uniqueness
Theorems could not be applied. The key point of the proof they gave was just that the image of a certain cycle
without exits was again a cycle, possibly with more edges than the original one. The hypotheses in this case were 
less general than the ones in the following generalization of the Cuntz-Krieger Uniqueness Theorem.

\begin{theor}\label{variante} Let  $E$ be an arbitrary graph, $A$ a graded $K$-algebra and  $\pi: L_K(E)\to A$ a
ring homomorphism  with $\pi(v) \neq 0$, for every vertex $v\in E^0$, which maps each cycle without exits to a
non-nilpotent homogeneous element of nonzero degree. Then $\pi$ is injective.
\end{theor}
\begin{proof}
Note that the kernel of $\pi$ is an algebra ideal of $\L$ which does not contain vertices. If $Ker (\pi)$ is
nonzero, by Corollary \ref{corzuru} it contains a nonzero element $p(c,c^*)$, where $p$ is a polynomial and
$c$ is a cycle without exits. By the hypothesis $\pi(c)=h\ne 0$ is a homogeneous element of degree $r\neq 0$,
thus $0=\pi(p(c,c^*))=p(h,h^*)$. Since $h$ is not nilpotent, then the coefficients of the polynomial
$p(c,c^*)$ are all zero, a contradiction.
\end{proof}

Finally we show how to use Proposition \ref{zuru} to simplify the proof on the characterization of simple
Leavitt path algebras (see \cite[Theorem 3.1]{AA3}).

\begin{corol} Let $E$ be an arbitrary graph. Then $\L$ is simple if and only if $E$ satisfies Condition
{\rm (L)} and the only hereditary and saturated subsets of $E^0$ are the trivial ones.
\end{corol}
\begin{proof}
If $\L$ is simple, both conditions on the graph $E$ are  proved in \cite[Theorem 3.1]{AA3}. For the converse
take into account that Condition (L) implies that any nonzero element in $\L$ is linked to a vertex (see
Proposition \ref{VerticeEnTodoIdeal}). Thus, there is a vertex in any nonzero ideal $I$ of $\L$. But on the
other hand $\emptyset\ne I\cap E^0$ is hereditary and saturated (\cite[Lemma 2.3]{AA3}), therefore it
coincides with $E^0$ and so $I=L_K(E)$.
\end{proof}

%%%%%%%%%%%%%%%%%%%%%%%%%%%%%%%%%%%%%%%%%%%%%%%%%%%%%%%%%%%%%%%%%%%%%%%%%%%%%%%%%%%%%%%%%%%%%%%%%%%%%%%%%%%
%%%%%%%%%%%%%%%%%%%%%%%%%%%%%%%%%%%%%%%%%%%%%%%%%%%%%%%%%%%%%%%%%%%%%%%%%%%%%%%%%%%%%%%%%%%%%%%%%%%%%%%%%%%
\section{Minimal left ideals}

Minimal left ideals are the building pieces of the socle of a semiprime ring. Clearly enough, in order to be
able to compute $\soc{L_K(E)}$, it would be wise to collect as much information as possible on the structure
of these ideals. Hence, the aim of this section is to find necessary and sufficient conditions so that a
principal left ideal is minimal.

As a first step, we will find necessary and sufficient conditions on a vertex so that the left ideal it
generates turns out to be minimal. We recall some notions introduced and results proved in \cite{AMMS} for
row-finite graphs, which will be also useful in the context of arbitrary graphs.

We say that a \emph{path} $\mu$ \emph{contains no bifurcations} if the set $\mu^0\setminus\{r(\mu)\}$ contains
no bifurcations, that is, if none of the vertices of the path $\mu$, except perhaps $r(\mu)$, is a
bifurcation.

The following two results are valid verbatim for arbitrary graphs because the use of relation (4) in their
proofs is limited to the case of vertices $v$ without bifurcations (and therefore finite-emitters).

\begin{lem}\label{nobifurcationsutov}{\rm \cite[Lemma 2.2]{AMMS}} Let $E$ be an arbitrary graph and let $u,v$
be in $E^0$, with $v\in T(u)$. If there is only one path joining $u$ with $v$ and it does not contain
bifurcations, then $L_K(E)u\cong L_K(E)v$ as left $L_K(E)$-modules.
\end{lem}

Note that the following proposition assumes that $u$ is a finite-emitter.

\begin{prop}\label{prior}{\rm \cite[Proposition 2.3]{AMMS}} Let $E$ be an arbitrary graph and $u$ a regular
vertex with $s^{-1}(u)=\{f_1,\ldots,f_n\}$. Then $L_K(E)u=\bigoplus_{i=1}^n L_K(E)f_if_i^*$. Furthermore, if
$r(f_i)\ne r(f_j)$ for $i\ne j$ and $v_i:=r(f_i)$, then $L_K(E)u\cong\bigoplus_{i=1}^n L_K(E)v_i$.
\end{prop}

The next result, however, requires a slight adaptation from its row-finite analog.

\begin{lem}\label{priorinfinite} Let $E$ be an arbitrary graph and let $u\in E^0$ be an infinite emitter. Then
$\bigoplus_{i=1}^\infty L_K(E)f_if_i^*\subsetneq L_K(E)u$, where
$s^{-1}(u)=\{f_i\}_{i\in {\mathbb N}}$.  In particular, $L_K(E)u$ is not a
minimal left ideal.
\end{lem}
\begin{proof}
The inclusion $\bigoplus_{i=1}^\infty L_K(E)f_if_i^*\subseteq L_K(E)u$ is clear. Suppose that $u\in
\bigoplus_{i=1}^\infty L_K(E)f_if_i^*$ and write $u=\sum_j \alpha_j g_jg_j^*$, where $g_j\in s^{-1}(u)$. Since
$s^{-1}(u)$ is infinite, there exists $f\in s^{-1}(u)$ such that $f\neq g_j$ for all $j$. Then, $f=uf=\sum_j
\alpha_j g_jg_j^*f=0$, a contradiction.
\end{proof}

Recall that a left ideal $I$ of an algebra $A$ is said to be \emph{minimal} if it is nonzero and the only
left ideals of $A$ that it contains are $0$ and $I$. From the results above we get an immediate consequence.

\begin{corol} \label{lemmaA} Let $E$ be an arbitrary graph and $w\in E^0$. If $T(w)$ contains some
bifurcation, then the left ideal $L_K(E)w$ is not minimal.
\end{corol}

Thus we have found a first necessary condition for the minimality of the left ideal generated by a vertex. But
as in the row-finite case, there is a second condition, introduced in \cite{AMMS}. The proof given there holds
also in our more general setting.

\begin{prop}\label{closedpathminimal}{\rm \cite[Proposition 2.5]{AMMS}} Let $E$ be an arbitrary graph.
If there is some closed path based at $u\in E^0$, then $L_K(E)u$ is not a minimal left ideal.
\end{prop}

Thus using this proposition and Corollary \ref{lemmaA} we conclude:

\begin{prop}\label{sufficient} Let $E$ be an arbitrary graph. Let $u$ be a vertex of the graph $E$ and
suppose that the left ideal $L_K(E)u$ is minimal. Then $u\in P_l(E)$.
\end{prop}

As we shall prove in what follows, this necessary condition turns out to be also sufficient. Following the
reasoning given in \cite[Proposition 2.7]{AMMS} but using Corollary \ref{lemmaA} and Propositions
\ref{nodegenerada} and \ref{sufficient} instead, we have:

\begin{prop}\label{corner} Let $E$ be an arbitrary graph. For any $u\in E^0$, the left ideal $L_K(E)u$ is
minimal if and only if $uL_K(E)u=K u\ \cong K$.
\end{prop}

\begin{rem}\label{sink}
{\rm For any sink $u$, trivially $uL_K(E)u=K u\cong K$, and therefore the left ideal $L_K(E)u$ is minimal.
Also, if $w$ is a vertex connected to a sink $u$ by a path without bifurcations, then we have that $L_K(E)w$
is a  minimal left ideal because $L_K(E)w\cong L_K(E)u$ by Lemma \ref{nobifurcationsutov}.}
\end{rem}

Our task now is to show that the converse implication of Proposition \ref{sufficient} holds too. The
proof of this fact strongly differs from that given in the row-finite setting, and this is so precisely because we
lack the direct limit construction in which a great part of the proof for the row-finite case is based on. Our
new approach is more combinatorial.

Before proceeding with this task, we need to establish several preliminary results.

\begin{prop}\label{pro1} Let $E$ be an arbitrary graph and $u\in P_l(E)$. Then every nonzero element of the
form $f_1 \cdots f_k g_1^* \cdots g_p^*$ with  $r(g_p^*)=u$ and $s(f_1), r(f_i), r(g_i^*)\in T(u)$, is either
the vertex $u$ or it can be written as $g_{k+1}^* \cdots g_p^*$, with $1<k<p$.
\end{prop}

\begin{proof} We proceed by induction on the number of real edges $k$. If we have $f_1g_1^* \cdots g_p^*$,
 since $r(f_1)=s(g_1^*)=r(g_1)$ then $f_1=g_1$, and therefore
$f_1g_1^* \cdots g_p^*=g_1g_1^* \cdots g_p^*=g_{2}^* \cdots g_p^*$, by (4). Suppose the result is valid for
$k-1$. Consider $f_1 \cdots f_kg_1^* \cdots g_p^*$; by the induction hypothesis $f_{2} \cdots f_kg_1^* \cdots
g_p^*=g_{k}^* \cdots g_p^*$ so that $f_1 \cdots f_kg_1^* \cdots g_p^*=f_1g_k^* \cdots g_p^*=g_{k+1}^* \cdots
g_p^*$.
\end{proof}

\medskip

\begin{lem}\label{lemmacorrected} Let $E$ be an arbitrary graph and let $\mu,\nu\in E^*$, with
$l(\mu),l(\nu)\geq 1$, $s(\mu)=s(\nu)$ and such that for every $u\in \mu^0\cup \nu^0$ there are neither
bifurcations nor cycles at $u$. Then, $\mu \nu^*\neq 0$ implies $\mu \nu^*=s(\mu)$.
\end{lem}

\begin{proof} We prove it by induction on $l(\mu)+l(\nu)$. The base case is for $l(\mu)+l(\nu)=2$.
In this case we have $\mu=f_1$ and $\nu=g_1$, with $s(f_1)=r(g_1^*)=s(g_1)$, and since we have no bifurcations
at $s(f_1)$, necessarily $f_1=g_1$ and moreover, by (4) we get $f_1g_1^*=s(f_1)$.

Let us suppose the result holds for the cases with $l(\mu)+l(\nu)<n$, and prove it  for
$l(\mu)+l(\nu)=n$. Write $\mu=f_1\dots f_r$ and $\nu=g_s\dots g_1$. Note that by the hypothesis we have $r,s\geq
1$. Now, since $\mu \nu^*\neq 0$, then
$$(\S) \quad f_2\dots f_rg_1^*\dots g_{s-1}^*\neq 0.$$ In this situation we have three possibilities:

If $r=1$, then again having no bifurcations at $s(f_1)$ implies that $f_1=g_s$, and by $(\S)$  we get that
$r(f_r)=s(g_1^*)$, that is, $s(g_1^*)=r(g_s)=s(g_s^*)=r(g_{s-1}^*)$. In other words, $g_{s-1}\dots g_1$ is a
closed path based at $r(f_1)$, and therefore there exists some cycle based at this same vertex, contradicting
our assumption. If $s=1$ we may proceed analogously. Finally, for the case $r,s>1$ we are allowed to apply the
induction hypothesis on $(\S)$ with the paths $\mu'=f_2\dots f_r$ and $\nu'=g_{s-1}\dots g_1$ which of course
verify that $s(\mu')=s(\nu')$. Thus, $\mu'(\nu')^*=s(\mu')$ and consequently
$\mu\nu^*=f_1s(\mu')g_s^*=f_1f_1^*=s(f_1)$, again by using (4) and the fact that there are no bifurcations at
$s(f_1)$.
\end{proof}

We would like to determine the different types of monomials we might encounter in $L_K(E)u$ when $u\in
P_l(E)$. In order to do this, we make the following definitions. Define $L$ to be set of edges $f\in E^1$ so
that $s(f),r(f)\in T(u)$. It is clear that the sources of edges in $E^1\setminus L$ are not in in $T(u)$
because there are not bifurcations at any vertex of $T(u)$. Therefore, for all paths $\mu=f_1\cdots f_k$ with
$f_i\in L$ and $\tau= g_1\cdots g_p$ with $g_j\not\in L$, we have $\mu\tau=0$ because $r(\mu)\in T(u)$ but
$s(\tau)\not\in T(u)$. Moreover, the graph $T=(T(u),L,r|_{T(u)},s|_{T(u)})$ is a \emph{line graph} (that is, a
graph without bifurcations which has only $u$ as a source).

\begin{prop}\label{types} Let $E$ be an arbitrary graph and $u\in P_l(E)$. The monomials generating
$L_K(E)u$ as a $K$-vector space are of the following types:
\begin{enumerate}
\item $u$.
\item $f_1\cdots f_k$, with $r(f_k)=u$ and $s(f_1)\ne u$.
\item $f_1\cdots f_kg_1^*\cdots g_p^*$ with $s(f_1)\ne u$, $r(g_p^*)=u$, $f_i\not\in L$ and $g_i\in L$.
\item $g_1^*\cdots g_p^*$ with $r(g_p^*)=u$, $s(g_1^*)\ne u$ and $g_i\in L$.
\end{enumerate}
\end{prop}
\begin{proof}
Consider first a monomial in real edges. Then it is necessarily of one the types (1) or (2). Take next  a
monomial of mixed type with real and ghost edges. Then it must be of the form $f_1\cdots f_kg_1^*\cdots g_p^*$
with $r(g_p^*)=u$. In case $s(f_1)=u$, applying Lemma \ref{lemmacorrected} we fall again in case (1). So we
can proceed supposing $s(f_1)\ne u$.  If for some $i$ we have $f_i\in L$, then  $i=k$ or
$f_{i+1},\ldots,f_k\in L$. Since $f_i\cdots f_kg_1^*\cdots g_p^*$ is nonzero and all the edges are in $L$ (and
also $r(g_p^*)=u$) we have  $f_i\cdots f_kg_1^*\cdots g_p^*=g_j^*\cdots g_p^*$ for some $j$ so that $1\le j\le
p$, by Proposition \ref{pro1}.

In a similar fashion we can proceed with the elements of type (4) to see that $g_i\in L$.
\end{proof}

We have now all the technical ingredients in hand to prove that the necessary condition for a principal left
ideal generated by a vertex to be minimal given in Proposition \ref{sufficient}, is also sufficient.

\begin{theor}\label{lolicandido} Let $E$ be an arbitrary graph and $u\in E^0$. Then $L_K(E)u$ is a minimal
left ideal if and only if $u\in P_l(E)$.
\end{theor}
\begin{proof}
Let $u\in E^0$ such that $L_K(E)u$ is minimal. Then $u\in P_l(E)$ by Proposition \ref{sufficient}.

Now we prove the converse. Take $u\in P_l(E)$ and $0\neq z\in L_K(E)u$. We will show that $u\in L_K(E)z$. By
Proposition \ref{types} we may write $z=z_1+z_2+z_3+z_4$, where $z_i$ is a linear combination of monomials of
type (i) in Proposition \ref{types}. We distinguish four cases.

\underline{Case 1:} $z_1\neq 0$. In this situation $z_1=ku$ for some $0\neq k\in K$, so
$u=k^{-1}uz_1=k^{-1}uz\in L_K(E)z$.

\underline{Case 2:} $z_1=0$ and $z_4\neq 0$. Let $kt_1^*\cdots t_l^*$ be a nonzero monomial in $z_4$  with $l$
minimal. For every monomial $f_1\cdots f_k$ of type (2) appearing in $z_2$ we have that $ut_l\cdots t_1
f_1\cdots f_k=0$ since $CP(u)=\emptyset$. Pick  a nonzero monomial $f_1\cdots f_kg_1^*\cdots g_p^*$ of type (3)
of $z_3$. Then $ut_l\cdots t_1f_1\cdots f_kg_1^*\cdots g_p^*=0$ because $f_i\not\in L$ and $t_i\in L$.
Moreover, if we consider $g_1^*\cdots g_p^*$, a monomial of type (4) appearing in $z_4$, different from
$t_1^*\cdots t_l^*$, we have that $ut_l\cdots t_1 g_1^*\cdots g_p^*=0$ because otherwise there would exist a
closed path based at $u$ (observe that since $l$ is minimal, we necessarily have $l<p$). This shows that
$u=k^{-1}ut_l\cdots t_1z\in L_K(E)z$.

\underline{Case 3:} $z_1, z_4=0$ and $z_2\neq 0$. Choose $kt_1\dots t_l$, a nonzero monomial  in $z_2$ with
$l$ maximal. Note that $t_l^*\cdots t_1^*z_2=ku$ because for a nonzero monomial $f_1\cdots f_r$ different from
$t_1\cdots t_l$ appearing in $z_2$ we have that $t_l^*\cdots t_1^*f_1\cdots f_r\neq 0$ would imply $l>r$ and
consequently $CP(u)\neq \emptyset$, a contradiction.

Now, choose $f_1\cdots f_r g_1^*\cdots g_s^*$, a  monomial of type (3) and consider the element $x=t_l^*\cdots
t_1^* f_1\cdots f_r g_1^*\cdots g_s^*$. Distinguish the following three situations. First, if $r<l$, then
$x=t_l^*\cdots t_{r+1}^*g_1^*\cdots g_s^*=0$ since $CP(u)=\emptyset$. Second, if $r=l$ then $x=ux=ug_1^*\cdots
g_s^*=0$  as $CP(u)=\emptyset$. Finally, if $r>l$ then $x=t_l^*\dots t_1^*  f_1\cdots f_l f_{l+1}\cdots  f_r
g_1^*\cdots g_s^*=uf_{l+1}\cdots  f_r g_1^*\cdots g_s^*$, and this would imply $f_{l+1}\in L$, a
contradiction.

This proves that $t_l^*\cdots t_1^*z_3=0$ so that  $u=k^{-1}t_l^*\cdots t_1^*z_2=k^{-1}t_l^*\cdots t_1^*z\in
L_K(E)z$.

\underline{Case 4:} $z_1, z_2, z_4=0$ and  $z_3\neq 0$.  Choose $kt_1\cdots t_lh_1^*\cdots h_m^*$, a  nonzero
monomial in $z_3$ with $l$ minimal. Now we have two possibilities:
\begin{enumerate}
\item[(i)] There is some summand $f_1\cdots f_r g_1^* \cdots g_s^*$
of $z_3$ such that $f_1\cdots f_r\neq t_1\cdots t_l$ . If $r=l$ then $t_l^*\cdots t_1^* f_1\cdots f_r g_1^*
\cdots g_s^*=0$, whereas if $r>l$, we would get that $s(f_{l+1})\in T(u)$ so that $f_{l+1}\in L$, which is
impossible.
\item[(ii)] $z_3$ is the monomial $k t_1\cdots t_lh_1^*\cdots h_m^*$.
\end{enumerate}
Hence, in any case, $t_l^*\cdots t_1^*z$ is an element which is under the conditions in Case 2, and therefore $u\in
L_K(E)t_l^*\cdots t_1^*z \subseteq L_K(E)z$.
\end{proof}

We close this section with the result that states that minimal left ideals are generated by line points. The proof of this  fact
follows the same sketch that the proof of \cite[Theorem 3.4]{AMMS}, now using Lemma \ref{CicloSinSalidas},
Proposition \ref{zuru} and Theorem \ref{lolicandido} instead of their row-finite analogs.

\begin{theor}\label{reduction} Let $E$ be an arbitrary graph and let $x$ be in $L_K(E)$ such that
$L_K(E)x$ is a minimal left ideal. Then, there exists a vertex $v\in P_l(E)$ such that $L_K(E)x$ is isomorphic
(as a left $L_K(E)$-module) to $L_K(E)v$.
\end{theor}

%%%%%%%%%%%%%%%%%%%%%%%%%%%%%%%%%%%%%%%%%%%%%%%%%%%%%%%%%%%%%%%%%%%%%%%%%%%%%%%%%%%%%%%%%%
%%%%%%%%%%%%%%%%%%%%%%%%%%%%%%%%%%%%%%%%%%%%%%%%%%%%%%%%%%%%%%%%%%%%%%%%%%%%%%%%%%%%%%%%%%
\section{The socle of a Leavitt path algebra}

Leavitt path algebras are semiprime by Proposition \ref{nodegenerada}. This implies, in particular, that
their left and right socles agree, which enables us to speak of the socle without distinguishing sides. However, it is
more convenient for us to work with left ideals, hence we will obtain the socle as the sum of all minimal left ideals. 

Recall that a homogeneous component of
the socle is the sum of all minimal left ideals which are isomorphic among themselves. Each homogeneous
component is also a (two-sided) ideal and the sum of all of them is the socle. Having characterized in the
previous section the minimal left ideals, we can apply these results to characterize $\soc\L$.  First of all
we would like to give a generating set of vertices of the socle as a two-sided ideal.

\begin{prop}\label{leftsocle} For an arbitrary graph $E$ we have
that $\sum_{u\in P_l(E)}L_K(E)u\subseteq \soc{L_K(E)}$. The reverse containment does not hold in general.
\end{prop}
\begin{proof} Use Theorem \ref{lolicandido} to show that for any $u\in P_l(E)$, the left ideal $L_K(E)u$
is minimal and therefore it is contained in the socle.

The reverse containment is not true in general as shows the example given in \cite[Proposition 4.1]{AMMS}.
\end{proof}

As in the case of a row-finite graph, the socle of a Leavitt path algebra $L_K(E)$ is generated as a two-sided
ideal by $P_l(E)$, the set of line points. To prove this, we can follow the steps in the proof of
\cite[Theorem 4.2]{AMMS} but using Propositions \ref{nodegenerada}, \ref{reduction} and \ref{leftsocle} rather
than their row-finite versions, jointly with the fact that the ideal generated by a subset $H$ of $E^0$ agrees
with the ideal generated by the hereditary saturated closure of $H$ (see the first assertion of \cite[Lemma
2.1]{APS}).

\begin{theor}\label{socleofLPA} Let $E$ be an arbitrary graph. Then $\soc{L_K(E)}=I(P_l(E))=I(H)$, where $H$
is the hereditary and saturated closure of $P_l(E)$.
\end{theor}

This result has an immediate but useful corollary.

\begin{corol}\label{nonzerosocle} For an arbitrary graph $E$, the Leavitt path algebra $L_K(E)$ has
nonzero socle if and only if $P_l(E)\neq \emptyset$.
\end{corol}

If $R$ is a ring, we let ${\mathbb M}_\infty(R)$ denote the ring of matrices of countable size over $R$ with
only a finite number of nonzero entries.

\begin{examples} {\rm By using these results, we can compute the socle of some Leavitt path algebras of
not necessarily row-finite graphs.
\begin{enumerate}
\item Consider the \emph{infinite edges graph} $E_\infty$ given
by $$\xymatrix{{\bullet}^v  \ar[r]^{(\infty)} & {\bullet}^w}$$ Then, by \cite[Lemma 1.1]{AA3},
$L_K(E_\infty)\cong {\mathbb M}_\infty(K)\vee K$, where the latter denotes the unitization of ${\mathbb
M}_\infty(K)$. Thus Theorem \ref{socleofLPA} gives another way to show that $\soc{{\mathbb M}_\infty(K)\vee
K}={\mathbb M}_\infty(K)$ via the Leavitt path algebra approach because $\soc{{\mathbb M}_\infty(K)\vee
K}=\soc{L_K(E_\infty)}=I(P_l(E_\infty))=I(\{w\})={\mathbb M}_\infty(K)$, where the last equality can be
obtained by using the isomorphism defined in \cite[Lemma 1.1]{AA3}.

\item Take the \emph{infinite clock graph} $C_\infty$ $$\xymatrix{ & {\bullet} & {\bullet} \\
 & {\bullet}^v \ar@{.>}[ul] \ar[u] \ar[ur] \ar[r] \ar[dr] \ar@{.>}[d]
\ar@{}[dl] _{(\infty)} & {\bullet} \\ &  & {\bullet}}$$ By \cite[Lemma 1.2]{AA3} we know that
$L_K(C_\infty)\cong \bigoplus_{i=1}^\infty {\mathbb M}_2(K)\oplus KI_{22}$, where $I_{22}$ is the element in
$\prod_{i=1}^\infty {\mathbb M}_2(K)$ given by $I_{22}=\prod_{i=1}^\infty E_{22}$, and $E_{22}$ is the
standard $(2,2)$-matrix  unit in ${\mathbb M}_2(K)$. Thus, using again Theorem \ref{socleofLPA} and the
isomorphism given in \cite[Lemma 1.2]{AA3} we get $\soc{\bigoplus_{i=1}^\infty {\mathbb M}_2(K)\oplus
KI_{22}}=\soc{L_K(C_\infty)}=I(P_l(C_\infty))=I(C_\infty^0\setminus\{v\})=\bigoplus_{i=1}^\infty {\mathbb
M}_2(K)$.
\end{enumerate}
}
\end{examples}

The next corollary is a generalization of \cite[Corollary 4.4]{AMMS}. Recall that $L_K(1,\infty)$ is  the Leavitt
path algebra of the \emph{infinite rose graph} $R_\infty$ $$\xymatrix{{\bullet} \ar@(ul,ur)^{(\infty)} }$$
considered in \cite[Examples 3.1 (ii)]{AA3}.

\begin{corol}\label{soclematrixLeavitt} For all
$m,n\in\mathbb N\cup\{\infty\}$, $\soc{{\mathbb M}_m(L_K(1,n))}=0$.
\end{corol}
\begin{proof}
When both $m$ and $n$ are finite, we just apply \cite[Corollary 4.4]{AMMS}. In other cases we consider the
graph $E_n^m$ given by $$
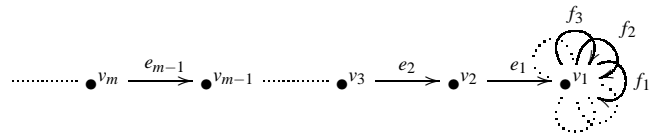
\xymatrix{\ar@{.}[r] & {\bullet}^{v_m} \ar [r] ^{e_{m-1}} & {\bullet}^{v_{m-1}}
\ar@{.}[r] & {\bullet}^{v_3} \ar [r] ^{e_2} &  {\bullet}^{v_{2}} \ar [r] ^{e_{1}} & {\bullet}^{v_1}
 \ar@(ur,dr) ^{f_1} \ar@(u,r) ^{f_2} \ar@(ul,ur) ^{f_3} \ar@{.} @(l,u) \ar@{.} @(dr,dl) \ar@{.} @(r,d)
 ^{\phantom {f_n}}& }$$ where if
$m=\infty$ then we have an infinite number of edges and vertices in the line, and if $n=\infty$ then we have an
infinite number of loops based at $v_1$. It is a tedious routine to check, with similar ideas to those of
\cite[Proposition 12]{AA2}, that ${\mathbb M}_m(L_K(1,n))\cong L_K(E_n^m)$. This
 graph satisfies that $P_l(E_n^m)=\emptyset$, for every $m,n\in\mathbb N\cup\{\infty\}$ so that Corollary
\ref{nonzerosocle} yields the result.
\end{proof}

We finish the paper by giving a structural characterization of the socle of a Leavitt path algebra $\L$ for an
arbitrary graph $E$.

If the socle is nonzero, then we know that $\soc\L=\oplus_\alpha C_\alpha$, where the $C_\a$ are the different
homogeneous components which are simple $K$-algebras (agreeing with their socles).

The structure theorem of simple algebras which coincide with their socles states that any such an algebra  is
isomorphic to an algebra $A=\F_{\Mp}(\M)$ (see \cite[IV, \S 8, p. 74]{Jacobson} for the definition). In the
framework of this theory $(\M,\M')$ is a  pair of dual vector spaces over a division $K$-algebra $\Delta$.
These vector spaces come from the minimal left ideal $\M=eA$ (so that $\M'=Ae$) and $\Delta$ is the division
$K$-algebra $\Delta=eAe$. Thus $\M$ is a left $\Delta$-vector space and $\M'$ a right $\Delta$-vector space.
Taking into account Proposition \ref{corner} we see that in our context $\Delta=K$ and the map $*\colon
\M\to\M'$ such that $ea\mapsto a^*e$ is an isomorphism of $K$-vector spaces. Hence $\dim(\M)=
\dim(\M')\le\aleph_0$. Then, applying \cite[IV, \S 15 Theorem 2, p. 89]{Jacobson} in the infinite-dimensional
case, each homogeneous component of $\soc\L$ is isomorphic to ${\mathbb M}_n(K)$, where $n\in\mathbb
N\cup\{\infty\}$.

Recall that a \emph{matricial algebra} is a finite direct product of full matrix algebras over $K$,  while a
\emph{locally matricial algebra} is  a direct limit of matricial algebras. Now, Litoff's Theorem \cite[IV, \S
15 Theorem 3, p. 90]{Jacobson} implies that each homogeneous component of the socle is locally matricial over
$K$ and so the socle itself is locally matricial over $K$. Thus we have proved the following

\begin{theor}\label{estructuradelzocalo}
For any arbitrary graph $E$ the socle of the Leavitt path algebra $L_K(E)$ is zero or a locally matricial algebra and
we have: $$\soc\L=\oplus_{n_i\in I}{\mathbb M}_{n_i}(K),$$ where $n_i\in{\mathbb N}\cup\{\infty\}$ and $I$ is a countable set.
\end{theor}

%%%%%%%%%%%%%%%%%%%%%%%%%%%%%%%%%%%%%%%%%%%%%%%%%%%%%%%%%%%%%%%%%%%%%%%%%%%%%%%%%%%%%%%%%%%%%%%%%%%%%%%%%%%%%%%
%%%%%%%%%%%%%%%%%%%%%%%%%%%%%%%%%%%%%%%%%%%%%%%%%%%%%%%%%%%%%%%%%%%%%%%%%%%%%%%%%%%%%%%%%%%%%%%%%%%%%%%%%%%%%%%
\section*{acknowledgments}

The authors wish to thank Pere Ara for valuable comments and correspondence. The first author was supported by
a Centre de Recerca Matem\`{a}tica Fellowship within the Research Programme ``Discrete and Continuous Methods
on Ring Theory". All authors were partially supported by the Spanish MEC and Fondos FEDER through projects
MTM2004-06580-C02-02 and MTM2007-60333, and by the Junta de Andaluc\'{\i}a and Fondos FEDER, jointly, through
projects FQM-336, FQM-1215 and FQM-2467. This work has also been supported by the Spanish Ministry of
Education and Science under project ``Ingenio Mathematica (i-math)'' No. CSD2006-00032 (Consolider-Ingenio
2010).

Part of this work was carried out during visits of the first author to the Universidad de M\'alaga and to the
Centre de Recerca Matem\`{a}tica. He thanks these host centers for their warm hospitality and support.

%%%%%%%%%%%%%%%%%%%%%%%%%%%%%%%%%%%%%%%%%%%%%%%%%%%%%%%%%%%%%%%%%%%%%%%%%%%%%%%%%%%%%%%%%%%%%%%%%%%%%%%%%%%%%%%
%%%%%%%%%%%%%%%%%%%%%%%%%%%%%%%%%%%%%%%%%%%%%%%%%%%%%%%%%%%%%%%%%%%%%%%%%%%%%%%%%%%%%%%%%%%%%%%%%%%%%%%%%%%%%%%


\begin{thebibliography}{99}
\bibitem{AA1} \textsc{G. Abrams, G. Aranda Pino}, The Leavitt path algebra of a graph, \emph{J. Algebra}
\textbf{293 (2)} (2005), 319--334.

\bibitem{AA2} \textsc{G. Abrams, G. Aranda Pino}, Purely infinite simple
Leavitt path algebras, \emph{J. Pure Appl. Algebra} \textbf{207 (3)} (2006), 553--563.

\bibitem{AA3} \textsc{G. Abrams, G. Aranda Pino}, The Leavitt path algebras
of arbitrary graphs, \emph{Houston J. Math.} (To appear.)

\bibitem{AAPS} \textsc{G. Abrams, G. Aranda Pino, F. Perera, M. Siles Molina}, Chain conditions for
Leavitt path algebras. (Submitted.)

\bibitem{AAS1} \textsc{G. Abrams, G. Aranda Pino, M. Siles Molina}, Finite-dimensional Leavitt path algebras,
\emph{J. Pure Appl. Algebra}   \textbf{209 (3)} (2007), 753-762.

\bibitem{AAS2} \textsc{G. Abrams, G. Aranda Pino, M. Siles Molina}, Locally finite Leavitt path algebras,
\emph{Israel J. Math.} (To appear.)

\bibitem{AGGP} \textsc{P. Ara, M.A. Gonz\'alez-Barroso, K.R. Goodearl, E. Pardo}, Fractional skew monoid rings,
\emph{J. Algebra} \textbf{278} (2004), 104--126.

\bibitem{AP} \textsc{P. Ara, E. Pardo}, Stable rank for graph algebras, \emph{Proc. Amer. Math. Soc.}
(To appear).

\bibitem{AMP} \textsc{P. Ara, M.A. Moreno, E. Pardo}, Nonstable K-Theory for graph algebras,
\emph{Algebra Represent. Theory}  \textbf{10 (2)} (2007), 157-178.

\bibitem{APS} \textsc{G. Aranda Pino, E. Pardo, M. Siles Molina}, Exchange Leavitt path algebras
and stable rank, \emph{J. Algebra} \textbf{305 (2)} (2006), 912--936.

\bibitem{Work} \textsc{G. Aranda Pino, F. Perera, M. Siles Molina, eds.}, \emph{Graph algebras: bridging
the gap between analysis and algebra}, ISBN: 978-84-9747-177-0, University of M\'{a}laga Press, M\'alaga,
Spain (2007).

\bibitem{AMMS} \textsc{G. Aranda Pino, D. Mart\'\i n Barquero, C. Mart\'\i n Gonz\'alez, M. Siles Molina},
 The socle of a Leavitt path algebra \emph{J. Pure Appl. Algebra} \textbf{212 (3)} (2008), 500-509.

\bibitem{ASgraded} \textsc{G. Aranda Pino, M. Siles Molina}, The Maximal
Graded Left Quotient Algebra of a Graded Algebra, \emph{Acta Math. Sinica, English Series}
 \textbf{22 (1)} (2006), 261--270.

\bibitem{BHRS} \textsc{T. Bates, J. H. Hong, I. Raeburn, W. Szyma\'{n}ski}, The ideal structure of the
C*-algebras of infinite graphs, \emph{Illinois J. Math.} \textbf{46 (4)} (2002), 1159--1176.

\bibitem{BPRS} \textsc{T. Bates, D. Pask, I. Raeburn, W. Szyma\'nski}, The C*-algebras of row-finite
graphs, \emph{New York J. Math.} \textbf{6} (2000), 307--324.

\bibitem{C} \textsc{J. Cuntz}, Simple C*-algebras generated by isometries, \emph{Comm. Math. Physics}
\textbf{57} (1977), 173--185.

\bibitem{CK} \textsc{J. Cuntz, W. Krieger}, A class of C*-algebras and topological Markov chains,
\emph{Invent. Math.} \textbf{63} (1981), 25--40.

\bibitem{DT} \textsc{D. Drinen, M. Tomforde}, The C*-algebras of arbitrary graphs,
\emph{Rocky Mountain J. Math} \textbf{35 (1)} (2005), 105--135.

\bibitem{FLR} \textsc{N. Fowler, M. Laca and I. Raeburn}, The $C^*$-algebras
of infinite graphs, \emph{Proc. Amer. Math. Soc.} \textbf{8} (2000), 2319-2327.

\bibitem{GarciSimon}\textsc{ J. L. Garc\' \i a, J.J.  Sim\'on}, Morita equivalence for idempotent rings.
\emph{J. Pure Appl. Algebra} {\textbf 76} {(1991)}, 39--56.

\bibitem{GS} \textsc{M. G\'omez Lozano, M. Siles Molina}, Quotient rings and Fountain-Gould left orders
by the local approach, \emph{Acta Math. Hungar.} \textbf{97} (2002), 287--301.

\bibitem{Goodearl} \textsc{K. R. Goodearl}, Leavitt path algebras and direct limits (To appear.)  arxiv.org/abs/0712.2554.

\bibitem{Jacobson} \textsc{N. Jacobson}, Structure of Rings, American Mathematical Society Colloquium
Publications, Vol. 37. 1968.

\bibitem{KPR} \textsc{A. Kumjian, D. Pask, I. Raeburn}, Cuntz-Krieger algebras of directed graphs,
\emph{Pacific J. Math.} \textbf{184 (1)} (1998), 161--174.

\bibitem{Ky} \textsc{S. Kyuno}, Equivalence of module categories. \emph{Math. J. Okayama Univ.}
\textbf{28} (1974), 147--150.

\bibitem{Le} \textsc{W.G. Leavitt}, The module type of a ring, \emph{Trans. Amer. Math. Soc.} \textbf{103}
(1962), 113--130.

\bibitem{R} \textsc{I. Raeburn},
\emph{Graph algebras}, CBMS Regional Conference Series in Mathematics, \textbf{103}, Amer. Math. Soc.,
Providence, (2005).

\bibitem{RS} \textsc{I. Raeburn, W. Szyma\'nski}, Cuntz-Krieger algebras of infinite graphs and matrices,
\emph{Trans. Amer. Math. Soc.} \textbf{356 (1)} (2004), 39--59.

\bibitem{SilesAQLPA} \textsc{M. Siles Molina},
Algebras of quotients of Leavitt path algebra, \emph{J. Algebra} (To appear).
DOI:10.1016/j.jalgebra.2007.09.017.

\bibitem{Tomforde} \textsc{M. Tomforde}, Uniqueness theorems and ideal structure for Leavitt path algebras,
 \emph{J. Algebra} \textbf{318 (1)} (2007), 270--299.


\end{thebibliography}
\end{document}